\documentclass[12pt,reqno]{amsart}
\usepackage{amsmath,amsthm,amsfonts,amssymb,amscd}
\usepackage{pxfonts}   

\numberwithin{equation}{section}

\newcommand{\Jac}{\operatorname{Jac}}

 \newcommand{\eps}{\varepsilon}

\newcommand{\eqm}{\overset{\circ}{=}}

\newcommand{\s} {\sigma}

\newcommand{\ti}{\pitchfork}
\newcommand{\U}{\mathcal U}
\newcommand{\V}{\mathcal V}

\def \cD {{\mathcal D}}

\def \cH {{\mathcal H}}

\def \cP {{\mathcal P}}

\def \cU {{\mathcal U}}

\def \cW {{\mathcal W}}

\newcommand \W {\cW}

\newtheorem{maintheorema}{Theorem}

\newtheorem{theorem}{Theorem}[section]

\newtheorem{corollary}[theorem]{Corollary}

\newtheorem{proposition}[theorem]{Proposition}

\newtheorem{lemma}[theorem]{Lemma}

\newtheorem{definition}[theorem]{Definition}

\newtheorem{obs}[theorem]{Remark}

\theoremstyle{remark}

\begin{document}

\thanks{}

\author{F. Rodriguez Hertz}
\address{IMERL-Facultad de Ingenier\'\i a\\ Universidad de la
Rep\'ublica\\ CC 30 Montevideo, Uruguay.}
\email{frhertz@fing.edu.uy}
\urladdr{http://www.fing.edu.uy/$\sim$frhertz}
\author{M. A. Rodriguez Hertz}
\address{IMERL-Facultad de Ingenier\'\i a\\ Universidad de la
Rep\'ublica\\ CC 30 Montevideo, Uruguay.}
\email{jana@fing.edu.uy}
\urladdr{http://www.fing.edu.uy/$\sim$jana}
\author{A. Tahzibi}
\address{Departamento de Matem\'atica,
  ICMC-USP S\~{a}o Carlos, Caixa Postal 668, 13560-970 S\~{a}o
  Carlos-SP, Brazil.}
\urladdr{http://www.icmc.sc.usp.br/$\sim$tahzibi}
\email{tahzibi@icmc.sc.usp.br}
\author{R. Ures}
\address{IMERL-Facultad de Ingenier\'\i a\\ Universidad de la
Rep\'ublica\\ CC 30 Montevideo, Uruguay.}
\email{ures@fing.edu.uy}
\urladdr{http://www.fing.edu.uy/$\sim$ures}

\thanks{This work was done in Universidad de la Rep\'{u}blica- Uruguay and ICMC-USP, S\~{a}o Carlos, Brazil.
 A. Tahzibi would like to thank FAPESP for financial support and Universidad de la Rep\'{u}blica for warm hospitality and
 financial support. F. Rodriguez Hertz and R. Ures also acknowledge warm hospitality of ICMC-USP
 and financial support of CNPq and FAPESP respectively.
 F. Rodriguez Hertz, M. A. Rodriguez Hertz, and Ra\'{u}l Ures were partially supported by
 PDT 54/18 and 63/204 grants and A. Tahzibi was also partially supported by CNPq (Projeto Universal)}

\date{}
\keywords{}

\subjclass{Primary: 37D25. Secondary: 37D30, 37D35.}

\renewcommand{\subjclassname}{\textup{2000} Mathematics Subject Classification}

\date{}

\setcounter{tocdepth}{2}

\title[A Criterion For Ergodicity]{A Criterion For Ergodicity of Non-uniformly hyperbolic Diffeomorphisms}

\begin{abstract}
In this work we exhibit a new criteria for ergodicity of
diffeomorphisms involving conditions on Lyapunov exponents and
general position of some invariant manifolds. On one hand we derive
uniqueness of SRB-measures for transitive surface diffeomorphisms.
On the other hand, using recent results on the existence of blenders
we give a positive answer, in the $C^1$ topology, to a conjecture of
Pugh-Shub in the context of partially hyperbolic conservative
diffeomorphisms with two dimensional center bundle.
\end{abstract}

\maketitle


\section{Introduction}

In this note we announce certain criteria to prove ergodicity, and
their application in two different settings: stable ergodicity and
uniqueness of SRB-measures.

\subsection{Stable ergodicity} Let $f: M
\rightarrow M$ be a volume preserving diffeomorphism of a compact
Riemannian manifold. We will denote by $m$ the probability Lebesgue
measure that is induced by a volume form $\omega.$ A challenging
problem in smooth ergodic theory is to prove the ergodicity of the
Lebesgue measure. Eberhard Hopf \cite{hopf} provided the first and
still the only argument of wide use to establish ergodicity and
proved the ergodicity  of the geodesic flow in the case of
negatively curved surfaces. Another interesting problem is the quest
of abundance and stability of ergodicity in the setting of volume
preserving diffeomorphisms. The celebrated Kolmogrov-Arnold-Moser
result, showing the presence of elliptic dynamics, is an obvious
obstruction to obtain density of ergodicity. Pugh and Shub proposed
a program to prove abundance (density) of stable ergodicity among
partially hyperbolic dynamical systems.

By Pugh-Shub conjecture we refer the density of stable ergodicity
among partially hyperbolic volume preserving diffeomorphisms.
Partial hyperbolicity is a weak form of hyperbolicity and partially
hyperbolic systems are far from elliptic dynamics.  A diffeomorphism
$f : M \rightarrow M$ is partially hyperbolic if $TM $ admits a non
trivial $Df$-invariant splitting $TM = E^s\oplus E^c \oplus E^u$,
such that all unit vectors $v^\s\in E^\s_x$ ($\s= s, c, u$) with
$x\in M$ verify:
$$\|D_xfv^s\| < \|D_xfv^c\| < \|D_xfv^u\| $$
for some suitable Riemannian metric. It is also required that
$\|Df|_{E^s}\| < 1$ and $\|Df^{-1}|_{E^u}\| < 1$. We denote the
subset of partially hyperbolic $C^{1+}-$diffeomorphisms of $M$ by
$\cP\cH (M)$ and $\cP\cH_m (M)$ represents the partially hyperbolic
diffeomorphisms preserving the volume Lebesgue measure. By a
$C^{1+}$ diffeomorphism we mean one with
H\"{o}lder continuous first order derivatives. \par%
 It is a known fact that there are foliations $\W^\s$ tangent to the distributions $E^\s$
for $\s=s,u$ (see for instance \cite{bp}).
 As a part of their program, Pugh and Shub conjectured that  transitivity of the pair of
 stable and unstable foliations, called accessibility, should be enough to prove ergodicity.
 That is, a diffeomorphism is accessible  if the unique nonempty set
 containing the stable and unstable manifold of every point is the whole
 $M$. By essential accessibility we mean that sets saturated by
 stable and unstable manifolds have full or null measure.
 Recently, Burns-Wilkinson \cite{burns-wilkinson2006} proved that essential accessibility and a
 bunching condition in the center direction imply ergodicity and consequently gave a positive
 answer to one of the Pugh-Shub's conjectures.
 Center bunching  is a technical condition that informally
 speaking  means quasi-conformality in the central direction.
 R. Hertz-R. Hertz-Ures \cite{rhrhu2006} proved the
 same result but for central direction of dimension one
 (where the center bunching condition is obviously verified)

 We remark that the first
examples of  stably ergodic diffeomorphisms which are not partially
hyperbolic were given in \cite{ta}.

 In the partially hyperbolic
context we apply our criterion to obtain stable ergodicity under
some conditions on central Lyapunov exponents and the existence of
some topological blenders. The novelty of this work is the use of
topological instruments to prove ergodicity.   We remark that an
antecedence of taking advantage of the interplay between nonuniform
and partial hyperbolicity is the following theorem of
Burns-Dolgopyat-Pesin \cite{burns-dolgopyat-pesin2002}.

\begin{theorem}[\cite{burns-dolgopyat-pesin2002}] \label{BDP}
Let $f\in\cP\cH_m (M)$. Assume that $f$ is essentially accessible
and has negative central exponents on an invariant set $A$ of
positive measure.
 Then $A$ has full measure and $f $ is stably ergodic.
\end{theorem}

 A natural question arises: what happens when the
exponents are different from zero but the signs may vary? In this
paper we consider this mixed case when the central dimension is two.

Finally, we show how to apply our results to prove Pugh-Shub
conjecture (in $C^1$ setting) for partially hyperbolic dynamical
systems with two dimensional central bundle.

\subsection{SRB-measures}

Another application of our criteria is  the case of certain
observable measures, in the sense that typical trajectories have
positive Lebesgue measure. Let $f$ be a diffeomorphism of a surface
$M$. We will say that an $f$-invariant ergodic measure $\mu$ is a
Sinai-Ruelle-Bowen measure (SRB-measure for short) if the largest
Lyapunov exponent is strictly positive and the Pesin´s formula holds
that is $\lambda^+(\mu)=h_\mu$ ($h_\mu$ is the entropy of the
measure $\mu$).

SRB-measures are important objects of study when Lebesgue measure is
not preserved. Since SRB-measures has a positive Lyapunov exponent,
unstable manifolds form  measurable partitions of their support.
Conditional measures along these partitions are absolutely
continuous with respect to the Lebesgue measure of the unstable
manifolds. See, for instance, the survey \cite{lsyoung2002} or
\cite{barreira-pesin} for  presentations of the subject.

An adaptation of our criterion (Theorem \ref{criterium}) to this
case allows us to show that transitive surface diffeomorphisms have
at most one SRB-measure (Theorem \ref{unique srb}). This contrasts
with known examples of transitive endomorphisms of surfaces having
intermingled basins of SRB-measures and with diffeomorphisms having
this property in greater dimensions (see \cite{kan1994} and
\cite[Ch. 11]{Bonatti-Diaz-Viana}).

\section{A Criterion for Ergodicity}

 Given a
hyperbolic periodic point $p$ we define the unstable manifold of
the orbit of $p$ as $W^u(O(p))=\bigcup_{k= 0}^{n(p)-1} W^u(f^k(p))$, similarly for
the stable manifold. Given a regular point, we define its Pesin's
stable manifold as usual, i.e.
$$
W^s_P(x)=\{y:\limsup_{n\to +\infty} \log
\frac{1}{n}d(f^n(x),f^n(y))<0\}$$ similarly for the unstable
manifold.

Given a hyperbolic periodic point $p$, let us define the following
sets:
\begin{eqnarray*}
B^s(p)=\{x:W^s_P(x)\ti W^u(O(p))\neq\emptyset\}\\
B^u(p)=\{x:W^u_P(x)\ti W^s(O(p))\neq\emptyset\}
\end{eqnarray*}
where $\ti$ means that the intersection is transversal.
$B^{\sigma}(p)$, $\sigma=s,u,$ is clearly $f-$invariant. By the
above definition any  $x \in B^s(p)$ has at least $\dim W^s(p)$
negative Lyapunov exponents and similarly for $x \in B^u(p).$
\begin{maintheorema}\label{criterium}
If $m(B^{\sigma}(p))>0$, for both $\sigma=s$ and $\sigma =u$ then,
$$
B(p):=B^u(p)\cap B^s(p)\eqm B^u(p)\eqm B^s(p),
$$
where the two last equalities are $m-$almost sure. Moreover,
$f|_{B(p)}$ is ergodic and non-uniformly hyperbolic and for $x\in
B(p)$, $\dim W^s_P(x)=\dim W^s(p)$ and $\dim W^u_P(x)=\dim
W^u(p)$.
\end{maintheorema}

Here, no partial hyperbolicity is required. Observe also that in the
above theorem we do not require all the exponents to be nonzero,
although some of them should be.  Also we do not require a priori
that $\dim W^s_P(x)=\dim W^s(p)$. Let us state an immediate
consequence of the $\lambda$-lemma.
\begin{lemma}
If $W^u(p)\ti W^s(q)\neq\emptyset$ then $B^u(p)\subset B^u(q)$ and
$B^s(q)\subset B^s(p)$.
\end{lemma}
So we have:
\begin{corollary}\label{criterio a usar}
If $W^u(p)\ti W^s(q)\neq\emptyset$, $m(B^u(p))>0$ and
$m(B^s(q))>0$, then $B(p)\eqm B(q)$. If in addition their union
has full measure then $f$ is ergodic and non-uniformly
hyperbolic.
\end{corollary}

We will apply the above criteria for the proof of the stable
ergodicity of Lebesgue measure for a special class of partially
hyperbolic systems. The novelty here is the use of blenders (robust
topological objects) to obtain stable ergodicity.

\begin{maintheorema}\label{t.2}
Let $f \in \cP\cH_m (M) $ satisfy the following properties:
 \begin{enumerate}
 \item $f$ satisfies the accessibility property,
 \item  there is a dominated splitting $E^c = E^- \oplus E^+$  of the central bundle into one-dimensional subbundles,
 \item $\int_M \lambda^- dm < 0 $ and $\int_M \lambda^+ dm > 0,$
 \item $f$ admits a $cs-$blender of stable dimension  $s+ 1$ and a $cu-$blender of unstable dimension  $u+1$
 and their periodic points are homoclinically related.
  \end{enumerate}
  then $f$ is stably ergodic.
\end{maintheorema}

Let us define $ \cP\cH_m
(M,2)$ the subset of partially hyperbolic conservative diffeomorphisms with two dimensional central bundle.
\begin{maintheorema}[Pugh-Shub Conjecture]\label{t.3}

There is a $C^1-$dense subset $\mathcal{D} \subset \cP\cH_m
(M,2) $ such that any $f \in\mathcal{D}$ is stably ergodic.
\end{maintheorema}

 \section{Blenders and Proof of theorem \ref{t.2}.}\label{section.t.2}

In this section we briefly state  the  definition and the main
properties of blenders that we need.
We refer the reader to
\cite{Bonatti-Diaz-Viana} for more details and references.

In order to give a definition of blender we adopt the operational
point of view of \cite{Bonatti-Diaz-Viana}.
\begin{definition}[\cite{Bonatti-Diaz-Viana}]\label{def.blenders}
Let $f:M\rightarrow M$ be a diffeomorphism and $Q$ be a hyperbolic
periodic point of unstable dimension $u+1$. We say that $f$ has a
$cu-$blender associated to $Q$ if there is a $C^1-$neighborhood
$\cU$ of $f$ and a $C^1-$open set $\cD$ of embeddings of an
$u$-dimensional disk $D^u$ into $M$, such that for every $g\in \cU$
every disk $D\in \cD$ intersects the closure of $W^s(Q_g)$, where
$Q_g$ is the analytic continuation of $Q$ for $g$. Define
$cs-$blender in an analogous way interchanging $u$ and $s$ and
unstable by stable.
\end{definition}

In fact, $(u+1)-$dimensional strips containing a disk of $\cD$ and
whose tangent spaces does not contain stable directions intersect
$W^s(Q_g)$. These strips are called {\em vertical strips} (for
details see \cite[Ch. 6.2]{Bonatti-Diaz-Viana}).

This property of blenders is the key point of our proof of
ergodicity. In fact by means of this geometrical property, we obtain
transversal intersection between stable and unstable manifolds even
if we do not have control on their size and shape.

\begin{obs}\label{obs.blenders}
Suppose that we have a $ cu$-blender of a partially hyperbolic
diffeomorphism associated to a periodic point $Q$ of unstable
dimension  $u+1$. The construction of blenders imply that there is
an open ball of size $\eps_0$ such that any disk of dimension $u+1$
through a point in this ball, whose tangent space is inside a cone
around $E^u\oplus E^+$ and containing a (large enough) strong
unstable disk must intersect $W^s(Q_g)$.
\end{obs}

\subsection{Proof of Theorem \ref{t.2}}
Let $f$ be a diffeomorphism as in Theorem \ref{t.2}. Observe that
conditions 2), 3) in Theorem \ref{t.2} are stable under
perturbations, hence there is a neighborhood $\U$ of $f$ in the
$C^1$ topology such that every $g$ in $\U$ satisfies these
conditions. Observe also that by condition 4), we may assume that
the boxes appearing in Remark \ref{obs.blenders} contain balls of
radius $\eps_0$ for every $g$ in this $\U$.

As in the proof of \cite[Theorem 2]{burns-dolgopyat-pesin2002}, we
have that since $f$ is accessible, there is a $C^1$ neighborhood
$\V$ of $f$ such that for every $g$ in this neighborhood,  every
$g$-orbit is $\eps_0$ dense. We shall see that every $C^{1+}$
diffeomorphism in $\V\cap\U$ is ergodic.

 Let $g$ be in $\U\cap\V$ and let us see that we are in the
hypothesis of Corollary \ref{criterio a usar}. So, we have $cs$- and
$cu$- blenders with boxes containing $\eps_0$ balls, $B_{cs}$ and
$B_{cu}$ respectively. Let $p$ be the periodic point associated to
$B_{cs}$ and $q$ be the periodic point associated to $B_{cu}$ (see
Definition \ref{def.blenders} and the discussion after it). We also
have that $W^s(O(p))$ intersects $W^u(O(q))$ transversally.

Let us see that for a.e. $x\in M$, if $\lambda^-(x)<0$ then $x\in
B^u(p)$. Since $B^u(p)$ is invariant it is enough to see that an
iterate $x'$ of $x$ is in $B^u(p)$. The $\eps_0$-density of a.e.
orbit implies that we can take $x'$ in $B_{cs}$. Moreover, since
$\lambda^-$ is invariant we have that $\lambda^-(x')<0$. So we have
that $W^u_P(x')$ contains a disk of dimension $(u+1)$, tangent to
$E^u\oplus E^-$ and containing a large strip in the direction of
$E^u$. Hence, by the main property of the blenders we have that
$W^u_P(x')$ intersects transversally $W^s(p)$ and hence $x'\in
B^u(p)$. Similarly we prove that for a.e. $x\in M$, if
$\lambda^+(x)>0$ then $x\in B^s(q)$.

Since the splitting $E^-\oplus E^+$ is dominated we have that
$\lambda^-(x)<\lambda^+(x)$ for a.e. $x\in M$, and hence either
$\lambda^-(x)<0$ or $\lambda^+(x)>0$. So we obtain that  either
$x\in B^u(p)$ or $x\in B^s(q)$ for a.e. $x\in M$. In other words
$B^u(p)\cup B^s(q)\eqm M$. On the other hand, since $\int_M
\lambda^-(x)dm<0$ we have that there is a set of positive measure
where $\lambda^-(x)<0$ for every $x$ in this set and hence
$m(B^u(p))>0$. Similarly $\int_M \lambda^+(x)dm>0$ implies that
$m(B^s(q))>0$. Then we are in the hypotheses of Corollary
\ref{criterio a usar} concluding that the system is ergodic and
non-uniformly hyperbolic.

 Finally observe, that any iterate of $g$ satisfies also the same
properties of $g$ and hence is also ergodic, so $g$ is Bernoulli.

\section{Pugh-Shub Conjecture}

In this section, we explain how to apply Theorem \ref{t.2} to prove
Pugh-Shub Conjecture in $C^1$ topology. Let $f\in\cP\cH_m$. We shall
perform a finite number of arbitrarily small $C^1$-perturbations,
among $C^{1+}$ volume preserving partially hyperbolic
diffeomorphisms, so that we obtain a stably ergodic diffeomorphism
arbitrarily close to $f$.

 Due to
\cite{dolgopyat-wilkinson2003}, we loose no generality in assuming
that $f$ is stably accessible. \par Let us call
$\lambda_f^-(x)\leq\lambda_f^+(x)$ the central Lyapunov exponents of
$x$ with respect to $f$, and recall that
$$\int_M (\lambda_f^-(x)+\lambda^+_f(x))\,dm (x)=\int_M \log \Jac (Df(x)|E^c_x)\,dm(x)$$
Observe that this amount depends continuously on $f$, due to
continuity of $E^c_x(f)$.  We have
\begin{theorem}[Baraviera-Bonatti
\cite{baraviera-bonatti2003}]\label{teo.baraviera.bonatti}
 Let
$f$ be a $C^1$ partially hyperbolic diffeomorphism, then there are
arbitrarily small $C^1$-perturbations $g$ of $f$ such that
$$\int_M\log \Jac(Dg(x)|E^c_x)dm(x) >\int_M\log \Jac(Dg(x)|E^c_x)dm(x). $$

\end{theorem}

This theorem allows us to choose a $C^1$-perturbation $g \in
\cP\cH_{m} (M)$ arbitrarily near to $f$ such that
\begin{equation}\label{suma.expoentes.positiva}
 \int_M ( \lambda^{-}_{g} (x) +  \lambda^{+}_{g} (x)) dm (x) > 0.
\end{equation}
or
\begin{equation}\label{suma.expoentes.negativa}
 \int_M ( \lambda^{-}_{g} (x) +  \lambda^{+}_{g} (x)) dm (x) < 0.
\end{equation}
We deal with the first case and the second case requires a similar
argument. Observe that the condition \ref{suma.expoentes.positiva}
is verified for any diffeomorphism $C^1-$close enough to $g$ . We
consider the two following cases:
\begin{enumerate}
\item $E^{c}$ does not admit a dominated splitting, \item There is
a dominated splitting of $E^c = E^{-} \oplus E^{+}.$
\end{enumerate}

If the first case occurs, we use a result of Bochi-Viana
\cite{bochi-viana2005} to obtain a new diffeomorphism $h$
$C^1$-close to $g$ such that $ \int_M  \lambda^- _h(x)dm(x)>0.$ This
implies that there exists a subset $A$ of $M$ with positive Lebesgue
measure such that $\lambda^-_h (x) > 0$ for all $x \in A.$
Consequently $\lambda^+ _h(x) > 0$ for $ x \in A$ and we are in the
setting of Theorem \ref{BDP}.\par

So, let us deal with the second case: $E^c$ admits a dominated
splitting. Note that (\ref{suma.expoentes.positiva}) implies that
$\int_M\lambda^+_g(x)dm(x)>0$. Theorem \ref{teo.baraviera.bonatti}
above implies that either $\int_M\lambda^-_g(x)dm(x)>0$ and Theorem
\ref{BDP} applies, or else $\int_M\lambda^-_g(x)dm(x)<0$. In this
last case, we want to produce a perturbation in such a way that we
are in the hypotheses of Theorem \ref{t.2}, that is, we want to find
$h\in\cP\cH_m(M)$ $C^1$-arbitrarily near $g$ admitting a
$cs$-blender of stable dimension $(s+1)$ and a $cu$-blender of
unstable dimension $(u-1)$ which are homoclinically related. We
begin by stating the following lemma:
\begin{lemma}
Let $f\in\cP\cH_m(M)$ be a stably accessible diffeomorphisms such
that $\dim E^c=2$. Then, either
\begin{enumerate}
\item $f$ is $C^1$-approximated by stably ergodic diffeomorphisms,
or \item $f$ is $C^1$-approximated by $g\in\cP\cH_m(M)$ having three
hyperbolic periodic points with stable dimension  $s$, $(s+1)$ and
$(s+2)$, respectively, where $s=\dim E^s$.
\end{enumerate}
\end{lemma}

To prove the above lemma we apply conservative versions of Ma\~{n}\'{e}'s
Ergodic Closing Lemma \cite{arnaud2001} and Frank's Lemma
\cite{bdp2003}.

To finish the proof we will need the following version of the
Connecting Lemma. A proof of more general results can be found in
\cite{arnaud2001} or \cite{boncrov2004}.

\begin{theorem}[Connecting Lemma]\label{connlemma}
Let $p, q$ be hyperbolic periodic points of a $C^{1+}$ conservative
transitive diffeomorphism $f$. Then, there exists a $C^{1+}$
conservative diffeomorphism $g$ $C^1-$close to $f$ such that
$W^s(p)\cap W^u(q)\neq \emptyset$.
\end{theorem}

We are almost done.  Since the diffeomorphism is transitive,
Connecting Lemma  implies, by making a perturbation, the existence
of cycles in the conditions of the following  proposition. This proposition is a conservative version of results in  \cite{bonatti-diaz2006}.

\begin{proposition}  \label{aproximateblender} Let $f$ be a $C^{1+}-$conservative
diffeomorphism having a co-index one cycle with real central
eigenvalues associated to saddles. Then $f$ can be approximated (in
the $C^1-$topology) by $C^{1+}-$conservative diffeomorphisms which
robustly ($C^1-$ topology) admit blenders.
\end{proposition}

Here, a co-index one cycle is a cycle where the difference between
the stable dimensions of the saddles is one.

After this we obtain the desired blenders and by applying the
Connecting Lemma again, we have that their periodic points are
homoclinically related. Finally, Theorem \ref{t.2} implies stable
ergodicity.

\section{Uniqueness of SRB-measures}

As we have said in the introduction we can also apply our criteria
to show uniqueness of SRB-measures of surface diffeomorphisms.

\begin{maintheorema}\label{unique srb}
Let $f$ be a $C^{1+}$ transitive diffeomorphism of a surface. Then,
$f$ has at most one SRB-measure.
\end{maintheorema}

As a corollary we also obtain that

\begin{corollary}
Let $f$ be a non-uniformly hyperbolic $C^{1+}$ {\sl conservative}
 transitive diffeomorphism of a surface. Then, $f$ is ergodic.
\end{corollary}

Let us say a few words about the proof of Theorem \ref{unique srb}.
First of all let us consider the easier case where $f$ satisfies the
Kupka-Smale (KS) condition. Suppose that $p_i, \, i=1,2$, are two
(typical) periodic points (as the ones obtained in \cite{katok1980})
associated to SRB-measures $\mu_i, \, i=1,2$. We can take
``rectangles" $R_i,\, i=1,2$, with sides in $W^\sigma(p_i), \,
i=1,2,\,\, \sigma=s,u$. Transitivity of $f$ implies that $f^k(R_1)$
intersects $R_2$ for some large enough $k$. If the rectangles $R_i$
have been chosen small enough the intersection between $f^k(R_1)$
and  $R_2$ implies that an iterate of $W^u (p_1)$ intersects
(transversally, since $f$ verifies KS condition) $W^s (p_2)$. Now,
analyzing the stable ``foliations" as in our criterion, the
equivalence of the measures $\mu_i$ to the Lebesgue measure along
unstable manifolds will imply that $\mu_1=\mu_2$.

In case $f$ does not verify the KS condition the strategy is to show
that many transversal intersections exist in spite of the possible
existence of tangencies. This is obtained by a subtle argument using
Sard´s Theorem.

\end{document}